
\input amstex
\documentstyle{amsppt}

\def\<{\hskip-2pt}
\def\({\hskip1pt(}
\def\){)\hskip1pt}
\def\>{\hskip1pt}
\def\ov{\overline}
\def\wh{\widehat}
\def\wt{\widetilde}
\def\pa{\partial}
\def\lgl{\longleftarrow}

\def\op{\operatorname}
\def\?{\hskip-1.2pt,\hskip-0.6pt}
\def\a{\alpha}
\def\be{\beta}
\def\ga{\gamma}

\def\th{\theta}

\def\la{\lambda}
\def\si{\sigma}
\def\ph{\varphi}
\def\om{\omega}
\def\Ga{\Gamma}
\def\De{\Delta}

\def\Si{\Sigma}

\def\La{\Lambda}

\def\Z{{\Bbb Z}}
\def\R{{\Bbb R}}
\def\C{{\Bbb C}}

\topmatter
\title Quantization and Coherent States over Lagrangian Submanifolds
\endtitle
\author Mikhail~V.~Karasev \endauthor
\address
Moscow Institute of Electronics and Mathematics,
3/12 B.~Vuzovsky per., \linebreak Moscow 109028, Russia;
e-mail: karasev\@aplmat.miemstu.msk.edu, root\@onti.miem.msk.su
\endaddress

\abstract
A membrane technique, in which the symplectic and Ricci
forms  are integrated over surfaces in a complexification of the phase
space, as well a ``creation" connection with zero curvature over lagrangian
submanifolds, is used to obtain a unified quantization
including a noncommutative algebra of functions, its representations, the
Dirac axioms, coherent integral transformations for solutions of spectral
and Cauchy problems, and trace formulas. \endabstract
\endtopmatter

\document
The results obtained in \cite{1} for the phase space $\R^{2n}$ and Gaussian
coherent states (see also \cite{2--5}) are developed in this work for
curved phase spaces and general coherent states.

\subhead \S1. Complexification
\endsubhead
Let $\frak X$ be a manifold endowed with a complex structure $J$. By
$(x,y)$ we denote points from the product $\frak X^\#=\frak X\times\frak
X$, and we equip this space with the complex structure $J\ominus J$ and
with the groupoid multiplication $(x,y)\circ (y,w)=(x,w)$. The set of units
(diagonal) in $\frak X^\#$ is identified with $\frak X\approx \op{diag}$;
the left and the right reductions $\frak X
\overset\pi\to\longleftarrow \frak X^\#
\overset\ov\pi\to\longrightarrow \frak X$ in this groupoid are given by the
projections $\pi(x,y)=x$, $\ov\pi(y,x)=x$. These mappings are,
respectively, a morphism and an antimorphism of complex structures, and so
$\frak X^\#$ is a {\it complex groupoid\/} corresponding to the manifold
$(\frak X, J)$, see  \cite{6, 7}. Let $\Pi(x)$,
$\ov\Pi(x)\subset \frak X^\#$ be fibers of $\pi$ and $\ov\pi$ over
$x\in\frak X$. Following physicists, we denote $(x,y)\equiv y|x$.

Note that one can identify the complexified tangent spaces ${}^\C T_x\frak
X$ with the tangent spaces $T_{x|x}\frak X^\#$ by using the natural
identifications $T_x\frak X\approx T_{x|x}\op{diag}$ and $i\cdot T_x\frak
X \approx(J\ominus J)T_{x|x}\op{diag}$. Then the eigenspaces of $J_x$
corresponding to the eigenvalues $+i$ and $-i$ coincide with the tangent
subspaces $T_{x|x}\Pi(x)$, $T_{x|x}\ov\Pi(x)$ in $T_{x|x}\frak X^\#$. Thus,
the groupoid $\frak X^\#$ can be regarded as a {\it complexification\/}
of the manifold $(\frak X,J)$; the fibers $\Pi(x)$ can be interpreted as
the integral leaves of the complex polarization $\Pi$ on $\frak X$
corresponding to the complex structure $J$, and the fibers $\ov\Pi(x)$ can
be interpreted as the leaves of the conjugate polarization $\ov\Pi$. The
groupoid inversion mapping $x|y\to y|x$ can be interpreted as an
involution in the complexification.

We use $z=\{z^k\}$ to denote local holomorphic coordinates on $\frak X$,
and by $\pa$ we denote the holomorphic differential, $\pa=\pa/\pa z$.
Thus,  $\Pi$ is generated by the vectors $\pa/\pa\ov z^k$.

\subhead \S2. Membrane amplitudes
\endsubhead
Throughout the sequel $\frak X$ is assumed connected and
simply connected. We call $\frak X$ a {\it symplec\-tic-K\"ahler manifold\/}
if it is equipped with a symplectic form $\om$ compatible with the complex
structure $J$, and is also endowed with a K\"ahlerian metric $g$. Let us
choose a lagrangian submanifold $\La\subset\frak X$ and suppose that
$$
\frac{1}{2\pi} \int_\Si \Big(\frac{\om}{\hbar}+\rho\Big)
-\frac{1}{2\pi}\oint_{\pa\Si}\nu_\La\in\Z,\qquad \forall\Si\subset\frak
X,\,\, \pa\Si\subset\La,
\tag 1
$$
where $\hbar>0$, $\rho=(i/2)\ov\pa \pa\op{ln}\op{det}g$ is one half  the
Ricci form on $\frak X$, and $\nu_\La$ is a {\it fundamental\/} $1$-form on
$\La$ such that $d\nu_\La=\rho\big|_\La$; see the details in \cite{7--9}. By
$\Si$ we denote {\it membranes\/} in $\frak X$, i.e., immersed oriented
two-dimensional surfaces with boundary.

We shall also be interested in several kinds of membranes in the
complexification $\frak X^\#$.

{\it Triangle.}
Connect points $\a,\be\in\La$, and $\be|\a\in\frak X^\#$  by paths
$\a\overset\Ga\to\lgl \be\lgl \be|\a\lgl\a$ going along $\Pi(\a)$,
$\ov\Pi(\be)$, and $\La$ successively. The last path along $\La$ will be
denoted by $\Ga=\Ga(\be|\a)$, and
$\text{\it triangle\/}_\Ga(\be|\a)$ will denote a membrane in $\frak X^\#$
with this circle boundary.

{\it Quadrangle.}
Connect points   $x\lgl x|y\lgl y\lgl y|x\lgl x$  along
$\Pi(x)$, $\ov\Pi(y)$, $\Pi(y)$, and $\ov\Pi(x)$ successively, and denote
by {\it quadrangle\/}$(y,x)$ a membrane in $\frak X^\#$ with this
boundary.

{\it Pentagon.}
Connect points   $\a\overset\Ga\to\lgl\be \lgl \be|x \lgl x
\lgl x|\a \lgl \a$ along $\Pi(\a)$, $\ov\Pi(x)$, $\Pi(x)$,
$\ov\Pi(\be)$, and $\La$. A membrane with this boundary is denoted by
$\text{\it pentagon\/}_\Ga(\be|x|\a)$.

{\it Hexagon.}
Paths $x\lgl x|w\lgl w\lgl w|y \lgl y \lgl y|x \lgl
x$  along $\Pi(x)$, $\ov\Pi(y)$, $\Pi(y)$, $\ov\Pi(w)$,
$\Pi(w)$, $\ov\Pi(x)$ generate a membrane $\text{\it hexagon\/}(w,y,x)$.

To each of the above membranes we relate a {\it membrane amplitude\/} as
follows:
$$
\align
a_\La(\be|\a) &= \exp\bigg\{i \int_{\text{\it triangle\/}_\Ga(\be|\a)}
\Big(\frac{\om}{\hbar}+\rho\Big) - i
\int_{\Ga(\be|\a)}\nu_\La\bigg\}\cdot
\De_\La(\a)^{1/2}\De_\La(\be)^{1/2},\\
 p(y,x) &= \exp\bigg\{i
\int_{\text{\it quadrangle\/}(y,x)}
\Big(\frac{\om}{\hbar}+\rho\Big)\bigg\},\\
k_\La(\be|x|\a) &= \exp\bigg\{i
\int_{\text{\it pentagon\/}_\Ga(\be|x|\a)} \Big(\frac{\om}{\hbar}+\rho\Big)
- i \int_{\Ga(\be|\a)}\nu_\La\bigg\}\cdot
\De_\La(\a)^{1/2}\De_\La(\be)^{1/2},\\
b(w,y,x) &= \exp\bigg\{i
\int_{\text{\it hexagon\/}(w,y,x)} \Big(\frac{\om}{\hbar}+\rho\Big)\bigg\}.
\endalign
$$
Here $\De_\La$ is a {\it modular function\/} on $\La$  \cite{8}, and
the form $\om/\hbar +\rho$ is assumed continued analytically on
$\frak X^\#$, possibly with singularities, so that all amplitudes $a_\La,
p, k_\La$, and $ b$ are globally defined and smooth. In view of
\thetag{1}, these amplitudes do not depend on the choice of membranes.
Note that \thetag{1} implies the condition
$1/2\pi[\om/\hbar+\rho]\in H^2(\frak X,\Z)$, well known in the
Kostant-Souriau  quantization.
We call $a_\La$ a {\it transition amplitude\/} and $k_\La$ a {\it coherent
amplitude\/} over $\La$; we also call $p$ and $b$ a {\it probability\/} and
a {\it holonomy\/} function over $\frak X$ \cite{8, 9}.

\subhead \S3. Geometric star product and the completeness axiom
\endsubhead
The probability and holonomy functions determine operators $\Cal
P$ and $\Cal B$ by the formulas
$$
\align
\Cal P f(x) &\overset{\text{def}}\to = \int_{\frak X} p(x,y)
f(y)\,dm(y),\\
\Cal B (f\otimes g) &\overset{\text{def}}\to = \int_{\frak X}\int_{\frak X}
b(w,y,x) f(y)g(x)\,dm(y)\,dm(x),
\endalign
$$
where $dm$ is a measure on $\frak X$.
The kernel $\op{Ker}\Cal P\subset L^2(\frak X,dm)$ is
nontrivial in general but the range of $\Cal B$ is orthogonal to
$\op{Ker} P$, so that the element
$$
f*g\overset{\text{def}}\to = \Cal P^{-1}\Cal B(f\otimes g),\qquad\forall
f,g\in C_0^\infty(\frak X)
\tag 2
$$
is well defined in the orthogonal complement $(\op{Ker}\Cal P)^\perp$. The
multiplication $*$ defined by \thetag{2} is noncommutative and
associative. The completion $\Cal L_{\frak X}$ of the space
$(\op{Ker}\Cal P)^\perp$ equipped with this multiplication and with the
inner product $(f,g)_{\Cal L_{\frak X}} \overset{\text{def}}\to= (\Cal
Pf,g)_{L^2(\frak X,dm)}$ is a Hilbert algebra; see details in \cite{6,
10}. This quantum algebra of functions on $\frak X$ has the unit $1$,
i.e., $1*f=f*1=f$, if the probability function $p$ and the measure $dm$
satisfy
$$
\int_{\frak X}p(x,y)\,dm(y)=1, \qquad \forall x\in\frak X.
\tag 3
$$
This is a {\it completeness\/} axiom.
Certain classes of homogeneous
K\"ahler manifolds automatically satisfy this axiom with respect to the
Liouville measure. In general the existence of a measure on $\frak X$ with
respect to which property \thetag{3} would be held is an open question.

In what follows we suppose that \thetag{3} holds.



There is a simple formula for the multiplication \thetag{2} in terms of
local complex coordinates. First, note that one can express the measure
$dm$ via local coordinates, $dm=M\,d\ov z\,dz$, and define the {\it dual
measure\/}  $dm' \overset{\text{def}}\to= M'\,d\ov z\,dz$ using the
density $M'=M^2/\op{det}g$. If there
is a K\"ahlerian metric $g'$ on $\frak X$ such that $M'=\op{det}g'$, then
$g '$ is referred to as a {\it dual metric\/}.

One can represent the symplectic form as  $\om=ig^\om_{k
s}\,d\ov z^k \wedge dz^s$, where $g^\om_{k
s}=\ov \pa_k \pa_s F$; here $F$ is a local potential. The corresponding
Liouville measure and ``symplectic" Laplace operator on $\frak X$ are
denoted by
$
dm^\om =\op{det} g^\om\cdot d\ov z\,dz$ and $\,{}^\om\De =
2(g^\om)^{-1}\,\ov\pa\,\pa.
$
Let us take a function $H\in C^\infty(\frak X)$ and define its {\it local
holomorphic realization\/} $\Cal H=\Cal H(\xi,z)$ by setting
$H= M^{\prime-1/2}\circ \Cal H(\pa F-\hbar\pa,z)(M^{\prime 1/2})$, where
$\pa\equiv \pa/\pa z$.
If $\Cal H(\xi,z)$ is a polynomial in the momenta $\xi$, then
$\Cal H(\pa F-\hbar\pa,z)$ is well defined as a Weyl-symmetrized operator.

\proclaim{Theorem 1}
For functions $H$ whose local holomorphic realizations $\Cal H$ are
polynomial in the momenta, the operator of left multiplication in the
algebra \thetag{2} is given by the formula
$$
H*=M^{\prime-1/2}\circ\Cal
H(\pa F-\hbar\pa,z)\circ M^{\prime 1/2}= H-i\hbar\, ad(H)_+ + \hbar^2
D(H)_{++}.
$$
Here $ad(H)_+$ is the $\pa$-part of the Hamilton field $ad(H)$, and
$D(H)_{++}$ is an operator regularly depending on $\hbar$ and annihilating
antiholomorphic functions.  \endproclaim

\subhead \S4. Evaluation of the star product. Holomorphic filtration
\endsubhead
Let $ad(H)$ denote the Hamilton field related to $H$ and $\om$. Let $\Cal
F^{(1)}(\frak X,\Pi)$ be the subspace of functions whose Hamilton flows
preserve the polarization, i.e., $H\in\Cal F^{(1)}\Leftrightarrow
[ad(H),\Pi]\subset\Pi$. By analogy, we define the space $\Cal F^{(2)}$ as
follows: $H\in\Cal F^{(2)}\Leftrightarrow [[ad(H),\Pi],\Pi]\subset\Pi$,
and so on. The sequence $\Cal F^{(1)}\subset \Cal F^{(2)}\subset\dots$
will be called a {\it holomorphic filtration\/}. The first space $\Cal
F^{(1)}$ is a Lie algebra with respect to the Poisson brackets.

Let us (locally) define   a $2$-tensor $((\Bbb
H^{sk}))$ by the relation  $\{\{H,f\},g\}=
\Bbb H^{sk}\pa_s f \pa_k g$ for any holomorphic functions $f$ and $g$.
This tensor is the $\ov\pa\pa$-part of the first variation of
$ad(H)$. Note that $H\in\Cal F^{(1)} \Leftrightarrow \Bbb H\equiv0$, and
$H\in\Cal F^{(2)}\Leftrightarrow\Bbb H$ is holomorphic. We consider three
special cases:

(I) $H\in\Cal F^{(1)}$; (II) $H\in\Cal F^{(2)}$, and the dual metric $g'$
on $\frak X$ is Ricci-flat; (III)  $H\in\Cal F^{(2)}$ and $\Bbb H$ is
nondegenerate, i.e., all local tensors $((\Bbb H^{ks}))$ are invertible.

In cases (I) and (II), let us define a {\it quantum deformation\/} as
follows \cite{7}:
$$H_\hbar \overset{\text{def}}\to= H- {\hbar} {}^\om
\De(H)/{4} - {i\hbar}\, ad(H)_+ (\op{ln}|{\Cal D m'}/{\Cal
D m^\om}|/2).
$$

In  case (III)  the inverse $\Bbb H^{-1}$ transforms under
changes of coordinates $z$ like a metric. Thus, there is the
``Laplace operator" ${}^{\Bbb H}\De$ with respect to this ``metric". Let
us denote
$$
\Bbb D(H) = \frac12(\op{det}\Bbb
H)^{-1/4}\circ{}^{\Bbb H}\De\circ(\op{det}\Bbb H)^{1/4},\quad
{}^{\Bbb H}\De \overset{\text{def}}\to= (\op{det}\Bbb
H)^{1/2}\pa_s\circ(\op{det}\Bbb H)^{-1/2}\Bbb H^{sk}\circ \pa_k, \tag 4
$$
and define  a secondary quantum deformation
$
H_{\hbar\hbar}\overset{\text{def}}\to= H_\hbar
+{\hbar^2}(M')^{-1/2} \Bbb D(H)(M^{\prime 1/2}).
$

\proclaim{Theorem 2}
{\rm(I)} If $H\in\Cal
F^{(1)}(\frak X,\Pi)$, then  $H_\hbar* = H_\hbar-i\hbar\, ad(H)_+$.

{\rm(II)} If $H\in\Cal F^{(2)}(\frak X,\Pi)$ and the dual metric $g'$ on
$\frak X$ is Ricci-flat, then
$
H_\hbar\, * = H_\hbar -i\hbar\, ad(H)_+ +{\hbar^2}(\op{det}
g')^{-1} \pa_s\circ (\op{det} g')\Bbb H^{sk}\circ\pa_k$.

{\rm(III)} If $H\in\Cal F^{(2)}(\frak X,\Pi)$ and the
$\ov\pa\pa$-variation of the field $ad(H)$ is nondegenerate, then
$
H_{\hbar\hbar}\,* = H_\hbar -i\hbar\, ad(H)_+
+{\hbar^2}(M')^{-1/2}\circ \Bbb D(H)\circ (M')^{1/2}$.
\endproclaim

The statement (I), as well as the following  Corollary, was obtained in
\cite{7}.

\proclaim{Corollary 1}
The quantum deformation is a homomorphism of the Lie
algebra\linebreak $\Cal F^{(1)}(\frak X, \Pi)$ into the associative algebra
\thetag{2}:
$
H_\hbar* G_\hbar -G_\hbar* H_\hbar = -i\hbar\{H,G\}_\hbar$, $\forall
H,G\in\Cal F^{(1)}$.
\endproclaim

\subhead \S5. Geometric inner product and quantization over lagrangian
submanifolds
\endsubhead
The transition and coherent amplitudes determine integral operators $\Cal
A_\La$ and $\Cal K_\La$ over $\La$:
$$
\Cal A_\La\ph(\be) \overset{\text{def}}\to= \int_\La
a_\La(\be|\a)\ph(\a)\,d\si(\a),\quad
(\Cal K_\La(x)\ph)(\be) \overset{\text{def}}\to=
\int_\La k_\La(\be|x|\a)\ph(\a)\,d\si(\a),
$$
where $d\si$ is a measure on $\La$. The kernel $\op{Ker}\Cal
A_\La\subset L^2(\La,d\si)$  is nontrivial in general, but the range of
$\Cal K_\La(x)$ $(\forall x\in\frak X)$ is orthogonal to $\op{Ker}\Cal
A_\La$, so
that the  operators
$$
\Cal C_\La(x) \overset{\text{def}}\to= \Cal A_\La^{-1}\Cal K_\La(x)
\tag 5
$$
are well defined in the orthogonal complement $(\op{Ker}\Cal A_\La)^\perp$.
The completion of  $(\op{Ker}\Cal
A_\La)^\perp$ with respect to the inner product
$(\ph_1,\ph_2)_\La \overset{\text{def}}\to= (\Cal
A_\La\ph_1,\ph_2)_{L^2(\La,d\si)}$ will be denoted by $\Cal L_\La$; see the
details in \cite{8} and \cite{9, 11}.

\proclaim{Theorem 3} Let $\La$ be a closed lagrangian submanifold in a
symplectic-K\"ahler manifold $\frak X$, and let conditions \thetag{1} and
\thetag{3} hold.
Then

{\rm(i)} the mapping $\Cal C_\La:\frak X\to\op{Hom}\Cal L_\La$ defined
in \thetag{5} is a quantization of $\frak X$ represented in the Hilbert
space $\Cal L_\La$ of functions on $\La$. This means that
\cite{7}:

-- $\Cal C_\La(x)$  is a one-dimensional orthogonal projection in
$\Cal L_\La$ for any $x\in \frak X$;

-- the first quantum correlator generated by $\Cal C_\La$ coincides with
the probability function over $\frak X$, i.e., $\op{tr}(\Cal C_\La(x)
\Cal C_\La(y)) = p(x,y)$;

-- the family $\Cal C_\La$ is smooth and complete, i.e.,
$\int_{\frak X}\Cal C_\La(x)\allowmathbreak dm(x)=I$.

{\rm(ii)} The dimensions are given by the
formulas $\op{dim}\Cal L_\La=(\op{dim}\Cal L_{\frak X})^{1/2}=\int_{\frak
X}\,dm$.

{\rm(iii)} The correspondence between functions on $\frak X$ and operators
over $\La$, defined by
$$
H\to\overset\vee\to H,\qquad  \overset\vee\to H\overset{\text{def}}\to=
\int_{\frak X} H(x)\, \Cal C_\La(x) \,dm(x), \tag 6
$$
possesses the following properties:
$
\overset\vee\to 1=I,\quad \overset\vee\to H\overset\vee\to G =
\overset\vee\to {H*G},\quad \overset\vee\to H^* = \overset\vee\to
{\ov H},
$
where the star product $*$ is defined in \thetag{2}. Thus, the mapping
$H\to
\overset\vee\to H$ gives a representation of the quantum algebra $\Cal
L_{\frak  X}$ in the Hilbert space $\Cal L_\La$.

{\rm(iv)} For functions $H$ whose local holomorphic realizations $\Cal H$
are polynomial in the momenta, the operator \thetag{6} is given by the
formula:
$$
\overset\vee\to H = q(\a)^{1/2}\circ\Cal H(\pa F(\a)-\hbar \wt\pa,
z(\a))\circ q(\a)^{-1/2}. \tag 7
$$
Here $q(\a)\overset\text{def}\to=\Cal D z(\a)/\Cal D\si(\a)$, $\a\in\La$,
$\{z=z(\a)\}$ is the local equation of the submanifold $\La$, $\wt\pa
\overset\text{def}\to= (\pa z/\pa\a)^{-1*} \pa/\pa\a$. The operators $\pa
F - \hbar \wt\pa$ and $z(\a)$ in \thetag{7} are Weyl-symmetrized.
\endproclaim

Formulas \thetag{6}, \thetag{7} represent the general construction of
quantum operators over a lagrangian submanifold. For functions $H$ of
the first or second holomorphic filtration one can evaluate the operator
$\overset\vee\to H$ explicitly as a first- or second-order differential
operator on $\La$.

\proclaim{Theorem 4}
For $H\in\Cal F^{(1)}(\frak X, \Pi)$ the operator ${\overset\vee \to
H}_\hbar$ on $\La$ defined by \thetag{6} has the form:
$$
{\overset\vee \to
H}_\hbar = H\big|_\La-i\hbar\Big(v(H)+\frac12 \op{div}^\si v(H)\Big).
\tag 8
$$
Here  $v(H)$ is
the projection of $ad(H)$ on $\La$ along  $\Pi$, and
$\op{div}^\si$ denotes the divergence with respect to  $d\si$.
The mapping $H\to {\overset\vee \to H}_\hbar$  satisfies the Dirac
axioms:
$$
[{\overset\vee \to
H}_\hbar, {\overset\vee \to
G}_\hbar] = -i\hbar \overset\vee \to
{\{H,G\}}_\hbar,\qquad
\overset\vee \to
{H}_\hbar^* = {\overset\vee \to
{\ov H}}_\hbar,\qquad \forall H,G\in\Cal F^{(1)}(\frak X,\Pi).
\tag 9
$$
\endproclaim

The quantization  \thetag{8} on  $\Cal F^{(1)}$ and the
implication \thetag{8} $\Rightarrow$ \thetag{9} were obtained in \cite{5,
8} independently from the ``associative" quantization  \thetag{6}.
The general formulas \thetag{6} and \thetag{7} are valid in  more
complicated cases. Let us consider, for instance,  case (III)  listed in
\S4.

\proclaim{Theorem 5}
Suppose that $H\in\Cal F^{(2)}(\frak X, \Pi)$ and the $\ov\pa\pa$-variation
of the Hamilton field $ad(H)$ is nondegenerate.
Then the operator ${\overset\vee \to
H}_{\hbar\hbar}$ on $\La$ defined by \thetag{6} has the form
$$
{\overset\vee \to
H}_{\hbar\hbar} = H\big|_\La - i\hbar \Big( v(H) + \frac12 \op{div}^\si
v(H)\Big) + \frac{\hbar^2}{2} q^{1/2}\circ \wt{\Bbb D}(H)\circ q^{-1/2}.
\tag 10
$$
Here the operator $\wt{\Bbb D}(H)$ is defined by Eqs.~\thetag{4} with
$\Bbb H$ replaced by $\wt{\Bbb H} = \Bbb H\big|_\La$ and $\pa$ by $\wt\pa$
{\rm(}where $\wt\pa$ and $q$ are introduced in Theorem~{\rm 3 (iv)}{\rm)}.
\endproclaim

The last summand in \thetag{10} is a second-order operator
on $\La\subset\frak X$. Resembling operators, but defined only in local
charts on $\La\subset\R^{2n}$, were met in  Maslov's theory of
semiclassical approximations. Operators in a sense similar
to \thetag{10} are also appeared in the Blattner--Kostant--Sternberg and
Atiyah--Hitchin approaches to geometric quantization. It is remarkable
that by using operators of type \thetag{10} one can describe irreducible
representations of certain physically important algebras with quadratic
commutation relations \cite{12}.

\subhead \S6. Connection with zero curvature over lagrangian submanifolds
\endsubhead
Let $\bold P:\frak X\to\op{Hom}L$ be an abstract quantization of the
symplectic-K\"ahler manifold $\frak X$ represented in a Hilbert space $L$
(i.e., the properties mentioned in Theorem~3~(i) hold; see also
\cite{7}). Then we obtain a correspondence between functions on $\frak X$
and operators, in $L$,
$$
H\to\wh H,\qquad \wh H\overset{\text{def}}\to= \int_{\frak X} H(x)\bold
P(x)\,dm(x) \tag 11
$$
with the following properties:
$\wh 1 = I,\quad \wh H\wh G = \wh{H*G},\quad \wh H^* = \wh{\ov H}$,
where the star product $*$ is defined in \thetag{2}. In  Berezin's
terminology, the function $H$ is called the  contravariant symbol of
the operator $\wh H$ (see references and discussion in \cite{6, 7, 10}).

Let us consider the following $1$-forms on $\frak X$ with values in the
space of functions over $\frak X$:
$
a_x^+(\cdot)\overset{\text{def}}\to= \hbar \pa_x\op{ln} p(x,\cdot)$,
$a_x^-(\cdot)\overset{\text{def}}\to=\hbar \ov\pa_x \op{ln} p(x,\cdot)$,
where $x\in\frak X$ and $p$ is the probability function. By applying to
$a^+$ and $ a^-$ the quantization mapping \thetag{11}, one obtains
$1$-forms $\wh a^+$ and $ \wh a^-$ on $\frak X$ with values in the space of
operators acting in $L$. In view of the identities $\hbar \pa\bold P = \wh
a^+\bold P$, $\wh a^-\bold P=0$, we call  $\wh a^+$ and $\wh a^-$
the {\it creation\/} and {\it annihilation\/} forms on $\frak X$ (see
in \cite{7}).

\proclaim{Theorem 6}
Suppose that $\bold P$ is a quantization of the symplectic-K\"ahler
manifold $(\frak X, J, \om,\allowmathbreak g)$, represented in a Hilbert
space $L$, and let $\La\subset \frak X$ be a lagrangian submanifold with
property \thetag{1}. Then the following operator-valued $1$-form
$$
\wh\th \overset{\text{def}}\to=-\frac{1}{\hbar}\wh a^+\big|_\La
+i\nu_\La\cdot I \tag 12
$$
on $\La$ determines a connection with zero curvature and with identity
global holonomy in the trivial $L$-bundle  over $\La$ {\rm(}i.e., in
$\La\times L${\rm)}. In formula \thetag{12}, $I$ is the identity operator in
$L$, and $\nu_\La$ is the fundamental $1$-form on $\La$.

If one fixes a point $\a_0\in\La$  and a ``vacuum" vector $e_{\a_0}\in L$
such that $\bold P(\a_0) = e_{\a_0} \otimes e_{\a_0}^*$, then the parallel
translation of $e_{\a_0}$ along paths $\a_0\to\a$ on $\La$ by means of the
connection \thetag{12} gives a smooth family of vectors in $L$:
$e_\a=\underset\longleftarrow\to{\op{Exp}}\Big\{-\int_{\a_0}^\a
\wh\th\Big\}e_{\a_0}$,
$\a\in\La$
with the following properties: $\bold P(\a)=e_\a\otimes e_\a^*$, $\,\,\wh
a_\a^- e_\a=0$.
If one defines
$$
u_\a\overset{\text{def}}\to= \De_\La(\a)^{1/2}e_\a,
\tag 13
$$
then $(u_\a,u_\be)_L = a_\La(\be|\a)$, where $a_\La$ is the transition
amplitude and $\De_\La$ is the modular function on $\La$.
\endproclaim

\subhead \S7. Coherent integral transformation
\endsubhead
Using the family of vectors \thetag{13}, we determine the
integrals
$$
U_\La(\ph) \overset{\text{def}}\to= \int_\La \ph(\a) u_\a\,d\si(\a),
\quad
\ph\in C_0^\infty(\La).
\tag 14
$$
\proclaim{Theorem 7}
The transformation \thetag{14} generates a unitary isomorphism
$U_\La:{\Cal L}_\La\to L$ intertwining quantizations \thetag{6} and
\thetag{11}, i.e., $\wh H\circ U_\La= U_\La\circ \overset\vee\to H$.
\endproclaim

Following \cite{11}, we call vectors $u_\a\in L$ {\it geometric coherent
states\/} over $\La$ represented in the Hilbert space $L$ (or, in short,
$\La$-{\it coherent states\/} in $L$). The mapping $U_\La$ will be called a
{\it coherent integral transformation\/} over $\La$.

\proclaim{Corollary 2}
If $\ph$ is an eigenfunction of the operator $\overset\vee\to H$ in
$\Cal L_\La$, then $U_\La(\ph)$
is the eigenvector of the operator $\wh H$ in $L$ related to the same
eigenvalue. \endproclaim

This statement is the basis for constructing useful formulas as well as
for solving  equations in mathematical physics. If the phase space
$(\frak X, J, \om, g)$, the submanifold $\La\subset \frak X$ and the
quantization $\bold P:\frak X\to\op{Hom}(L)$ are taken in an appropriate
way, integral \thetag{14} leads to many well-known integral
representations of special functions.

     From the other hand the quantization procedure itself on the level of
operators, star products, etc., can be expressed in terms of
transformations $U_\La$ \thetag{14} and operators $\overset\vee\to H=
\overset\vee\to H_\La$ \thetag{6}. For instance, one has the following
formulas $\wh H=U_{\op{diag}}(H)$, $H*=\overset\vee\to H_{\op{diag}}$, where
$\op{diag}$ denotes  the lagrangian diagonal in $\frak X^\# = \frak
X\times\frak X$  (note that $\frak X^\#$ is, in a natural way, a
symplectic-K\"ahler groupoid); see in \cite{5, 7}.

The coherent transformation
\thetag{14} is also very effective for the theory of semiclassical
approximation (what was the basic observation and the starting point in
\cite{1}; see also \cite{2--5, 9, 13--17}). In the semiclassical approach
the number $\hbar$ is a parameter tending to zero over a certain subset
$\Cal R\subset (0,\infty)$, $0\in[\Cal R]$. In the general case, the metric
on $\frak X$ can depend on this parameter $g=g_{(\hbar)}$ as a regular series
in $\hbar\to 0$, and $\om$ is supposed to be the K\"ahler form for the
metric $g_{(0)}$. The rule \thetag{1} and axiom \thetag{3} are assumed for
each $\hbar\in\Cal R$. The measures $dm=dm_\hbar$ on $\frak X$ and $d\si =
d\si_\hbar$ on $\La$ must be represented as follows: $dm_\hbar =
(2\pi\hbar)^{-\op{dim}\frak X/2}(dm^\om + O(\hbar))$,
$d\si_\hbar=(4\pi\hbar)^{-\op{dim}\frak X/4}(d\si^0 + O(\hbar))$.

\proclaim{Theorem 8}
Under the above assumptions the following asymptotic expansions hold as
$\hbar\to0$.
For the star product \thetag{2}:
$
f*g = fg - i\hbar \,ad(f)_+ g + O(\hbar^2)$,
$f*g - g*f = -i\hbar\{f,g\} + O(\hbar^2)$.
For the quantum operators \thetag{6} over lagrangian submanifolds:
$
\overset\vee\to H= H\big|_\La - i\hbar(\,v(H)
+\op{div}^{\si^0}\! v(H)/2\,) + O(\hbar^2)$;
and, moreover, if $H|_\La\equiv\la_0=\op{const}$ and  $d\si^0$
 is invariant with respect to the Hamilton flow $\ga_H^\tau$, then
$\overset\vee\to H = \la_0-i\hbar d/d\tau + O(\hbar^2)$. In the latter
case the vector $u=U_\La(1)$ is an approximate eigenvector
{\rm(}quasimode{\rm)} for the operator $\wh H$ corresponding   to an
eigenvalue $\la_0+O(\hbar^2)$, i.e., $\wh H u =\la_0 u +O(\hbar^2)$,
$\|u\|^2 = \int_\La d\si^0 + O(\hbar)$.
\endproclaim

\subhead \S8.
Cauchy problem and trace formula
\endsubhead
The coherent transformation \thetag{14} can be applied not only to the
eigenvalue problem, but to the Cauchy problem as well \cite{2, 5, 16, 17}.
Let us take a real smooth function $H$ on $\frak X$. The graph
$(\ga_H^t)=\{(\ga^t_H(x),x)\mid x\in\frak X\}$ is a lagrangian submanifold
in $\frak X^\#=\frak X\times\frak X$. Let us transport the dual metric
$dm'$ from $\frak X$ onto the graph $(\ga_H^t)$ by means of the projection
$\frak X^\#\overset\ov\pi\to\longrightarrow \frak X$. The quantum operator
\thetag{6} related to the graph $(\ga_H^t)\subset\frak X^\#$ and to the
function $H^\#=H\otimes 1$ on $\frak X^\#$ can be represented as
follows:
$\overset\vee\to{H}^\# =  H -i\hbar\,ad(H)_+ +\hbar^2 D^t(H)_{++}$.
The latter operator $D^t(H)_{++}$ vanishes if $H\in\Cal
F^{(1)}(\frak X,\Pi)$, and it is a second order differential operator (in
$x$) if $H\in\Cal F^{(2)}(\frak X,\Pi)$; see Theorems~4,5. We denote by
$f^t(x)$ the solution of the Cauchy problem over the manifold $\frak X$:
$$
i\pa f^t/\pa t = \hbar D^t(H)_{++} f^t,\qquad f^t\big|_{t=0}=f.
\tag 15
$$
Then the solution of the Cauchy problem in the Hilbert space $L$ is
expressed as follows.

\proclaim{Theorem 9}
Let the operator $\wh H$ \thetag{11} be selfadjoint in the Hilbert space
$L$. Then for any $f\in C_0^\infty(\frak X)$ one has the formula
$
\exp\{-{it}\wh H/{\hbar}\}\wh f = U_{\text{\it
graph\/}(\ga_H^t)}(f^t) =  U_{\text{\it
graph\/}(\ga_H^t)}(f) + O(\hbar)$,
where the coherent integral transformation related to the $\text{\it
graph\/}(\ga_H^t)\subset \frak X\times\frak X$ is defined by \thetag{14},
and $f^t$ is defined by \thetag{15}. If $H\in\Cal F^{(1)}(\frak X,
\Pi)$, then $f^t\equiv f$.  \endproclaim

\proclaim{Corollary 3 (Trace formula)}
$$
\multline
\op{tr}\Big(e^{-it\wh H/\hbar}\wh f\Big) = \int_\frak X
\exp\Big\{i\int_{\Si^t(x)} \Big(\frac{\om}{\hbar}+\rho\Big) -
i\int_0^t \nu^\tau(x)\,d\tau\\
 - \frac{it}{\hbar}H(x)\Big\} \De^t(x)^{1/2}
f^t(x)\,dm(x).
\endmultline
$$
Here $\Si^t(x)$ is a triangle membrane with the boundary $x\longleftarrow
x|\ga^t(x)\longleftarrow \ga^t(x)\longleftarrow x$ {\rm(}see \S2, the
first path is the trajectory of the Hamilton flow{\rm)},  $\nu^t$ is a
fundamental form along the trajectory:
$$
\nu^t = \frac{i}{4}\frac{\pa}{\pa t} \op{ln} \frac{\Cal D z^t/\Cal
Dz}{\Cal D \ov z^t/\Cal D\ov z} -\frac12 \Im [(\pa\op{ln} \op{det}g)(z^t)
\dot z^t],\quad z^t\equiv z(\ga^t(x)),\quad z\equiv z(x),
$$
and $\De^t$ is a modular function along the trajectory:
$$
\De^t = |{\Cal Dz^t}/{\Cal D z}|
({\op{det}g(\ga^t(x))}/{\op{det}g(x)})^{1/2}.
$$
In particular, if the flow $\ga_H^t$ preserves the K\"ahlerian structure
on $\frak X$, then
$$
\op{tr}\Big(e^{-it\wh H/\hbar}\wh f\Big) =
\int_\frak X \exp\Big\{i\int_{\Si^t(x)}
\Big(\frac{\om}{\hbar}+\rho\Big) - i\int_0^t \nu^\tau(x)\,d\tau -
\frac{it}{\hbar}H(x)\Big\}  f(x)\,dm(x).
$$
\endproclaim

\enddocument